\documentclass[a4paper]{article}
\usepackage{amsmath,amsfonts,amssymb}
\usepackage{graphicx, epsfig}

\newtheorem{theorem}{Theorem}[section]

\newtheorem{prop}[theorem]{Proposition}

\newtheorem{rem}[theorem]{Remark}

\newenvironment{proof}[1][Proof]{\textbf{#1.} }
{\hfill\rule{0.5em}{0.5em}\medskip}
\newenvironment{proof*}[1][Proof]{\textbf{#1.} }{}

\def\epsilon{\varepsilon}

\begin{document}

\title{
Compact bordered Riemannian surfaces as vibrating membranes:
an estimate \`a la Hersch-Yang-Yau-Fraser-Schoen
\\}

\author{Alexandre Gabard}
\maketitle

\newbox\abstract
\setbox\abstract\vtop{\hsize 9cm \noindent

\footnotesize \noindent\textsc{Abstract.} We try to present an
estimate relating the first Dirichlet and Neumann eigenvalues
of a compact bordered Riemannian surface.}

\centerline{\hbox{\copy\abstract}}

\section{Introduction}\label{sec1}

In a recent paper Fraser-Schoen \cite{Fraser-Schoen_2011} took
advantage of Ahlfors' conformal representation of compact
bordered Riemann(ian) surfaces $\Sigma$ over the disc to
obtain a Steklov eigenvalue estimate
in terms of their topological invariants (the genus $p\ge 0$
and the number of contours $r\ge 1$, i.e. boundary
components). We do not need to recall here the rich history
antedating the result of Ahlfors 1950~\cite{Ahlfors_1950}
(presented already in
Spring
1948
 at Harvard\footnote{Compare, Nehari,
 1950, Trans. AMS, p.\,258.}), except for saying that
 this (pre)history is
(surprisingly?)
confined to the
{\it schlicht} case (i.e., $p=0$) which involves primarily a
contribution of Riemann\footnote{A Riemann's Nachlass worked
out by H. Weber, compare the historical interrogations raised
by Bieberbach in \cite{Bieberbach_1925}.} and subsequently
Schottky 1877, Bieberbach 1925, Grunsky 1937--41\footnote{For
accurate references
we refer the interested reader to the
bibliography in \cite{Gabard_2006}, which is by far not
exhaustive, e.g., a serious omission is R. Courant, {\it
Conformal mapping of multiply connected domains}, Duke Math.
J. 5 (1939), 814--823.}. In the present note we have attempted
to use the same method (as Fraser-Schoen) to get a similar
estimate for the classical vibrating membrane problem (e.g.,
Poisson 1829, Helmoltz 1862, Clebsch 1862, Lord Rayleigh
(=J.\,W. Strutt) 1894--96, Weyl 1911, Courant 1918,
Faber--Krahn 1923--24, etc.)\footnote{Precise references are
given in Kuttler--Sigillito, {\it Eigenvalues of the Laplacian
in two dimensions}, SIAM Review 26 (1984), 163--193.}:
\begin{equation}\label{vibrating_membrane}
- \Delta u = \lambda u \,,
\end{equation}
where $\Delta$ denotes the Laplacian (of Beltrami 1867)
attached to the Riemannian metric. As
the nature of the question
seems to impose it one must not focalize on the {\it fixed
membrane} (under Dirichlet boundary condition $u=0$ on
$\partial \Sigma$) nor on the {\it free membrane} problem
(under the Neumann boundary condition ${\partial u \over
\partial n}=0$ on $\partial \Sigma$, where $n$ is the normal to the boundary $\partial \Sigma$),
but rather more consider both problems in
some natural symbiosis (suggested by the Pythagorean geometry
of the sphere $x_1^2+x_2^2+x_3^2=1$). Then
the sought for estimate becomes very straightforward (indeed
completely parallel to Hersch 1970 \cite{Hersch_1970}). To
picture out the right historical perspective as a commutative
diagram, recall that Yang-Yau 1980
\cite[Prop.,\,p.\,58]{Yang-Yau_1980} generalized the first
estimate of Hersch \cite[Inequality~(1),
p.\,1645]{Hersch_1970} for Riemannian metrics on the sphere to
arbitrary closed (oriented) surfaces, whereas the present note
tries to achieve the same
goal regarding the second estimate of Hersch
\cite[Inequality~(2), p.\,1646]{Hersch_1970} involving
bordered surfaces topologically equivalent to the disc. Hoping
that the understanding of the (newcomer) author is
trustful, the key trick seems to use as ``isoperimetric''
model not the flat round disc but rather the (north)
hemisphere of the (unit) sphere, which ``sounds'' better.

{\small {\bf Numerical justification:} Indeed, comparing the
quantity $\lambda_1 A$ (where $\lambda_1$ is the first
Dirichlet eigenvalue, and $A=\mbox{area}$) we get for the disc
$j^2 \pi\approx 5.783 \pi$
(where $j\approx 2.4048255576$ is the first
positive zero of the
Bessel function $J_0$), while for the hemisphere we have
$2\cdot 2\pi\approx 4\pi$
which
has
a gravest fundamental tone (than the planar
disc). In contradistinction for the free
membrane problem the quantity $\mu_1 A$ ($\mu_1$=first nonzero
Neumann eigenvalue) has now to be maximized for a ``good
sounding''! We find for the disc $\mu_1 A= p^2 \pi\approx
3.390 \cdot \pi$, where $p\approx 1.8411837813$ is the first
positive zero of the
Bessel function $J_1'$, while for the hemisphere we have
$2\cdot 2\pi$ which is larger (hence ``better'').

}

\section{An inequality extending the one of Hersch as a
bordered avatar of the one by Yang-Yau}

\begin{prop} Let $\Sigma=\Sigma_{p,r}$ be a compact bordered
Riemannian orientable surface of genus $p$ with $r$ contours
of total area $A$. Denote by $\lambda_1$ the first Dirichlet
eigenvalue for the problem {\rm (\ref{vibrating_membrane})},
and by $\mu_1\le \mu_2$ the first two non-zero Neumann
eigenvalues. Assume the existence of a conformal mapping
$f\colon \Sigma \to D^2=\{ z\in {\Bbb C} : |z|\le 1\}$ to the
disc of degree $d$. Then we have the inequality
\begin{equation}\label{Hersch-Yang-Yau-inequality}
\bigl({1\over \lambda_1}+{1\over \mu_1}+{1\over \mu_2}\bigr)
{1\over A} \ge {1 \over d } {3 \over 4 \pi} \,.
\end{equation}
\end{prop}

\begin{rem} {\rm Ahlfors
\cite[\S 4, pp.\,122--133]{Ahlfors_1950} showed that for such
a Riemann surface there is always a holomorphic branched
covering to the disc $D^2$ of degree $\le r+2p$, whereas the
present author modestly improved the degree bound as being
$\le r+p$ (compare \cite{Gabard_2006}).
Hence the degree $d$ involved in inequality
(\ref{Hersch-Yang-Yau-inequality}) can
be taken as $r+p$. Of course for some particularized Riemann
surfaces one can hope to be more economical.
}
\end{rem}

\begin{rem} {\rm In the case where  the topology is simple
$\Sigma \approx D^2$ then by the Riemann mapping theorem we
may choose $d=1$ and inequality
(\ref{Hersch-Yang-Yau-inequality})
turns into an equality for the (unit) hemisphere $H=S^2 \cap
\{ x_3\ge 0\}$ as in this case $\lambda_1=\mu_1=\mu_2=2$.
(This is of course already observed in  Hersch
\cite{Hersch_1970}.)
}
\end{rem}

\begin{proof}
We merely
have to follow the idea of conformal transplantation of
P\'olya-Szeg\"o (1951), as elaborated
subsequently by Hersch 1970 \cite{Hersch_1970} and Yang-Yau
1980~\cite{Yang-Yau_1980}, conjointly with the variational
characterization of eigenvalues (Poincar\'e 1890, Rayleigh
1894, Fischer 1905, Ritz 1908, Courant 1920, P\'olya-Schiffer
1954, Hersch 1961 \cite{Hersch_1961})\footnote{Accurate
references as on p.\,100 of C. Bandle, {\it Isoperimetric
Inequalities and Applications}, Pitman, 1980.}. So let
$f\colon \Sigma \to D^2$ be our conformal mapping. As it will
be soon apparent it is more convenient to work with the
(north) hemisphere (instead of the flat disc)
$$
H=S^2\cap\{x_3\ge 0\} \quad\mbox{of}\quad
S^2=\{(x_1,x_2,x_3)\in{\Bbb R}^3 : x_1^2+x_2^2+x_3^2=1\}\,.
$$ The first (non-zero)
eigenvalues $\lambda_1$ and $\mu_1$ admits a variational
characterization as the absolute minimizers of the {\it
Rayleigh quotient}:
$$
R[u]={ \int_{\Sigma} |\nabla u|^2 dv \over \int_{\Sigma} u^2
dv}\,,
$$
where in the Neumann case orthogonality to the constant
functions (eigenfunctions for $0=\mu_0$) imposes the extra
side-condition $\int_{\Sigma} u dv=0$. Likewise Hersch
established in \cite{Hersch_1961} a variational
characterization for sums of reciprocals of eigenvalues. In
our situation this gives:
$$
{1\over \lambda_1}+{1\over \mu_1}+{1\over \mu_2}={\rm max}
(R[u_1]^{-1}+R[v_1]^{-1}+R[v_2]^{-1})\,,
$$
where $u_1$ satisfies the Dirichlet and $v_1,v_2$ the Neumann
boundary condition. The method is to pull-back (transplant via
$f$) the best functions on the target to get competitive trial
functions at the source. So on the hemisphere $H \subset {\Bbb
R}^3\ni(x_1,x_2,x_3)$ we consider the ambient coordinate
functions: $x_3$ verifying the Dirichlet condition, and $x_1,
x_2$ verifying the Neumann condition. The pull-backs $x_i\circ
f$ are eligible for the variational principle, since after
post-composing $f$ by a suitable conformal automorphism of the
hemisphere we may balance the {\it center of gravity}
$$
G=(\textstyle\int_\Sigma ( x_1\circ f)dv,\int_\Sigma (
x_2\circ f)dv)\in {\Bbb R}^2
$$
so as to make it coincide with the origin $(0,0)$. This
involves a topological argument initiated by Szeg\"o 1954
(later Weinberger 1956), which in our setting is (brilliantly)
exposed in Hersch 1970~\cite[Point 2.,
p.\,1646]{Hersch_1970}). We thus arrive at the inequality:
$$
{1\over \lambda_1}+{1\over \mu_1}+{1\over \mu_2}\ge
{\int_{\Sigma} (x_3\circ f)^2 dv
 \over \int_{\Sigma}|\nabla (x_3\circ f)|^2 dv}
+{\int_{\Sigma} (x_1\circ f)^2 dv
 \over \int_{\Sigma}|\nabla (x_1\circ f)|^2 dv}
 +{\int_{\Sigma} (x_2\circ f)^2 dv
 \over \int_{\Sigma}|\nabla (x_2\circ f)|^2 dv}\,.
$$
Each of the integrals  occurring in the denominators
$\int_{\Sigma} |\nabla (x_i\circ f)|^2 dv$ are equal to $d
\int_{H}|\nabla x_i|^2 dv={4\pi \over 3}$ (by conformal
invariance of the Dirichlet
integrand, compare Yang-Yau \cite[Lemma, p.\,59,
(ii)]{Yang-Yau_1980}). Adding up the numerators, we obtain, as
$\sum_{i=1}^3 (x_i\circ f)^2\equiv 1$ ($f$ taking values in
the unit sphere), finally $\int_{\Sigma} dv=A$. This complete
the proof of the proposed inequality
(\ref{Hersch-Yang-Yau-inequality}).
\end{proof}

One can also use merely the simple variational
characterization of the first eigenvalues to get first
$$
\mu_1 \int_\Sigma ( x_i \circ f)^2 dv \le \int_{\Sigma}|\nabla
(x_i\circ f)|^2 dv \quad (\mbox{for $i=1,2$})
$$
and likewise
$$
\lambda_1 \int_\Sigma ( x_3\circ f)^2 dv \le
\int_{\Sigma}|\nabla (x_3\circ f)|^2 dv \quad
\phantom{(\mbox{for $i=1,2$})}\,
$$
which added up (after multiplying by $\lambda_1$ the first two
inequalities and by $\mu_1$ the last one) lead to the
following
estimate involving only $\lambda_1$ and $\mu_1$:
\begin{equation}
\lambda_1 \mu_1 A \le d {4 \pi \over 3} (2 \lambda_1 +
\mu_1)\,.
\end{equation}
The latter inequality can of course also be deduced from
inequality (\ref{Hersch-Yang-Yau-inequality}) by using the
trivial inequation $\mu_1\le \mu_2$ (to eliminate $\mu_2$).

\medskip

{\small {\bf Acknowledgements.} The author  wishes to thank
several discussions with the musician Misha Gabard, the
Boul\'e family for financial support (i.e., sloshing
experiments with yoghurts in containers), as well as
conversations with Claude Weber, Michel Kervaire, Jean-Claude
Hausmann, Andr\'e Haefliger, Felice Ronga and
Johannes Huisman.}

{\small

}

{
\hspace{+5mm} 
{\footnotesize
\begin{minipage}[b]{0.6\linewidth} Alexandre
Gabard

Universit\'e de Gen\`eve

Section de Math\'ematiques

2-4 rue du Li\`evre, CP 64

CH-1211 Gen\`eve 4

Switzerland

alexandregabard@hotmail.com
\end{minipage}
\hspace{-25mm} }

\end{document}